\documentclass[10pt]{article}
\usepackage{amsmath,amssymb,amsthm,fullpage,graphicx,wrapfig}
\usepackage[a4paper, total={16cm, 23.7cm}]{geometry}

\newtheorem{theorem}{Theorem}[section]

\newtheorem{definition}[theorem]{Definition}

\numberwithin{equation}{section}

\begin{document}

\title{An introduction to\\stochastic processes associated\\with resistance forms\\and their scaling limits}
\author{D.~A.~Croydon\footnote{Department of Statistics, University of Warwick, Coventry, CV4
7AL, United Kingdom. \tt{d.a.croydon@warwick.ac.uk}}}
\maketitle

\begin{abstract} We introduce and summarise results from the recent paper `Scaling limits of stochastic processes associated with resistance forms' \cite{Croyres}, and also applications from `Time-changes of stochastic processes associated with resistance forms' \cite{CHKtime}, which was written jointly with T.\ Kumagai (Kyoto University) and B.~M.~Hambly (University of Oxford).
\end{abstract}

\section{Introduction}

The connections between electricity and probability are deep, and have provided many tools for understanding the behaviour of stochastic processes. In this note, we describe a new result in this direction from \cite{Croyres}, which states that if a sequence of spaces equipped with so-called `resistance metrics' and measures converge with respect to the Gromov-Hausdorff-vague topology, then the associated stochastic processes also converge. (In the non-compact case, the proof in \cite{Croyres} also requires a non-explosion condition.) All the relevant concepts will be introduced more carefully below, with the statement of the main result appearing as Theorem \ref{mainres}.

This result generalises previous work on trees, fractals, and various models of random graphs (apart from the background in \cite{Croyres}, see also \cite{ALWtree, kig1}). Moreover, it is useful in the study of time-changed processes, including Liouville Brownian motion, the Bouchaud trap model and the random conductance model, on such spaces \cite{CHKtime}. Some of these examples will be sketched in Section \ref{applications}. I further conjecture that the result will be applicable to the random walk on the incipient infinite cluster of critical bond percolation on the high-dimensional integer lattice (see Section \ref{IICsec}).

\section{Random walks on graphs and electrical networks}\label{rwsec}

Before introducing the definition of a resistance metric and associated stochastic process on a general space, it is helpful to recall the more elementary definition of effective resistance and the corresponding random walk on a graph. This is the purpose of the present section.

We start with the definition of a random walk on a weighted graph. In particular, let $G=(V,E)$ be a finite, connected graph, equipped with (strictly positive, symmetric) edge conductances $(c(x,y))_{\{x,y\}\in E}$. Let $\mu$ be a finite measure on $V$ of full-support. We then define the associated random walk $X$ to be the continuous time Markov chain with generator $\Delta$, as defined by:
\[(\Delta f)(x):=\frac{1}{\mu(\{x\})}\sum_{y:\:y\sim x}c(x,y)(f(y)-f(x)),\]
where the sum is over vertices $y$ connected to $x$ by an edge in $E$, i.e.\ this is the process that jumps from $x$ to $y$ with rate $c(x,y)/\mu(\{x\})$. Note that the transition probabilities of the jump chain of $X$ are given by \[P(x,y)=\frac{c(x,y)}{c(x)},\]
where $c(x):=\sum_{y:\:y\sim x}c(x,y)$, and so are completely determined by the conductances. The measure $\mu$ determines the time-scaling of the process. Common choices are to take $\mu(\{x\}):=c(x)$, which is the so-called \emph{constant speed random walk (CSRW)}, or $\mu(\{x\}):=1$, which is the \emph{the variable speed random walk (VSRW)}. As illustrated by the example presented in Section \ref{trapsec}, the latter two processes can have quite different behaviour if the conductances are inhomogeneous.

Suppose now we view $G$ as an electrical network with edges assigned conductances according to $(c(x,y))_{\{x,y\}\in E}$. If vertices in the network are held according to the potential $f(x)$, then the total \emph{electrical energy} dissipated in the network is given by $\mathcal{E}(f,f)$, where $\mathcal{E}$ is the quadratic form on $V$ given by
\[\mathcal{E}(f,g):=\frac12\sum_{x,y:x\sim y}c(x,y)\left(f(x)-f(y)\right)\left(g(x)-g(y)\right).\]
Moreover, regardless of the particular choice of $\mu$, $\mathcal{E}$ is a \emph{Dirichlet form} on $L^2(V,\mu)$, and can be written as
\[\mathcal{E}(f,g)=-\sum_{x\in V}(\Delta f)(x)g(x)\mu(\{x\}).\]
Using the classical correspondence between Dirichlet forms and reversible Markov processes, it follows that there is a one-to-one correspondence between the electrical energy $\mathcal{E}$ (viewed as a Dirichlet form on $L^2(V,\mu)$) and the random walk $X$. (For the definition of a Dirichlet form, and background on the connections between such objects and Markov processes, see \cite{FOT}.)

Suppose now that we wished to replace our network by a single resistor between two vertices $x$ and $y$. The resistance we should assign to this resistor to ensure that the same amount of current flows from $x$ to $y$ when voltages are applied to them as did in the original network is given by the \emph{effective resistance}, which can be computed by setting
\[R(x,y)^{-1}=\inf\left\{\mathcal{E}(f,f):\:f(x)=1,f(y)=0\right\}\]
for $x\neq y$, and $R(x,x)=0$. Although it is not immediate from the definition, it is possible to check that $R$ is a metric on $V$, e.g.\ \cite{Tet}, and characterises the edge conductances uniquely, e.g.\ \cite{Kigdendrite}. The latter observation is important, because it means that, given an effective resistance $R$ on a graph, one can reconstruct the corresponding electrical energy operator $\mathcal{E}$. Thus, if one is also given a measure $\mu$, then by viewing $\mathcal{E}$ as a Dirichlet form on $L^2(V,\mu)$ as in the previous paragraph we also recover the random walk $X$.

In summary, we have the following correspondences:
\begin{center}
\begin{tabular}{ccccc}
  Random walk $X$,   & $\leftrightarrow$ & Dirichlet form $\mathcal{E}$ &$\leftrightarrow$&Effective resistance $R$\\
  generator $\Delta$ &                   & on $L^2(V,\mu)$              &&
and measure $\mu$.
\end{tabular}
\end{center}

\section{Resistance metrics and forms}\label{ressec}

Building on the discussion of the previous section, it is now straightforward to introduce a resistance metric on a general space. After presenting the definition, we then explain the theory developed by Kigami in the context of analysis on low-dimensional fractals that links resistance metrics and stochastic processes (see \cite{kig1,Kig} for details).

\begin{definition}[{\cite[Definition 2.3.2]{kig1}}] Let $F$ be a set. A function $R:F\times F\rightarrow \mathbb{R}$ is a \emph{resistance metric} if, for every finite $V \subseteq F$, one can find a weighted (i.e.\ equipped with conductances) graph with vertex set $V$ for which $R|_{V\times V}$ is the associated effective resistance.
\end{definition}

\begin{figure}[h]
\begin{center}
\vspace{-0pt}
\scalebox{0.1}{\includegraphics{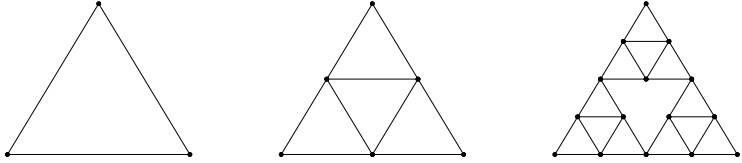}}
\caption{Level 0, 1, 2 approximations to the Sierpinski gasket.}\label{sgpic}
\end{center}
\end{figure}

As some first examples of resistance metrics, we have:
\begin{itemize}
  \item the effective resistance metric on a graph;
  \item the one-dimensional Euclidean metric $|x-y|$ on $\mathbb{R}$ (not true in higher dimensions), or fractional powers of this $|x-y|^{\alpha-1}$ for $\alpha\in(1,2]$ (see \cite[Chapter 16]{Kig});
  \item any `shortest path' metric on a tree-like metric space (see \cite{ALWtree, Kigdendrite});
  \item the resistance metric on the Sierpinski gasket, which can be constructed by setting, for `graph vertices' $x$, $y$ in the limiting fractal,
      \[R(x,y)=\lim_{n\rightarrow\infty}(3/5)^nR_n(x,y),\]
      where $R_n$ is the effective resistance on the level $n$ graph (see Figure \ref{sgpic}) considered with unit resistances along edges, and then using continuity to extend to whole space. Resistance metrics can similarly be defined on various classes of fractals, see \cite{kig1} for background.
\end{itemize}

Playing the role of the electrical energy in this general setting is the collection of resistance forms. We now state the definition of such objects. Whilst this is quite technical and we will not discuss the role of the various conditions in detail here, importantly it gives a route to connect the resistance metric with a stochastic process.

\begin{definition}[{\cite[Definition 2.3.1]{kig1}}] Let $F$ be a set. A pair $(\mathcal{E},\mathcal{F})$ is a resistance form on $X$ if it satisfies the following conditions:
\begin{description}
  \item[RF1] $\mathcal{F}$ is a linear subspace of the collection of functions $\{f:\:F\rightarrow \mathbb{R}\}$ containing constants, and $\mathcal{E}$ is a non-negative symmetric quadratic form on $\mathcal{F}$ such that $\mathcal{E}(f,f)=0$ if and only if $f$ is constant on ${F}$.
  \item[RF2] Let $\sim$ be the equivalence relation on $\mathcal{F}$ defined by saying $f\sim g$ if and only if $f-g$ is constant on $F$. Then $(\mathcal{F}/\sim,\mathcal{E})$ is a Hilbert space.
  \item[RF3] If $x\neq y$, then there exists an $f\in \mathcal{F}$ such that $f(x)\neq f(y)$.
  \item[RF4] For any $x,y\in F$,
  \[\sup\left\{\frac{|f(x)-f(y)|^2}{\mathcal{E}(f,f)}:\:f\in\mathcal{F},\:\mathcal{E}(f,f)>0\right\}<\infty.\]
  \item[RF5] If $\bar{f}:=(f\wedge 1)\vee 0$, then $f\in\mathcal{F}$ and $\mathcal{E}(\bar{f},\bar{f})\leq \mathcal{E}(f,f)$ for any $f\in\mathcal{F}$.
\end{description}
\end{definition}

The following theorem connects the notions of a resistance metric and a resistance form, and yields the stochastic process that will be of interest in the remainder of the article. In particular, it explains how the correspondences stated at the end of the previous section extend to the more general present setting. For simplicity of the statement, we restrict to the compact case. It is also possible to extend the result to locally compact spaces, though this requires a more careful treatment of the domain of the Dirichlet form.

\begin{theorem}[{\cite[Theorems 2.3.4, 2.3.6]{kig1}, \cite[Corollary 6.4 and Theorem 9.4]{Kig}}] (a) Let $F$ be a set. There is a one-to-one correspondence between resistance metrics and resistance forms on $F$. This is characterised by the relation:
\begin{equation}\label{reschar}
R(x,y)^{-1}=\inf\left\{\mathcal{E}(f,f):\:f(x)=1,f(y)=0\right\}
\end{equation}
for $x\neq y$, and $R(x,x)=0$.\\
(b) Suppose $(F,R)$ is compact resistance metric space, and $\mu$ is a finite Borel measure on $F$ of full support. Then the corresponding resistance form $(\mathcal{E},\mathcal{F})$ is a regular Dirichlet form on $L^2(F,\mu)$, and so naturally associated with a Hunt process $((X_t)_{t\geq 0},(P_x)_{x\in F})$.
\end{theorem}

As a first example of the connection between a resistance metric and a stochastic process (beyond the example of random walks on graphs already discussed), consider $F=[0,1]$, $R=\mbox{Euclidean}$, and $\mu$ be a finite Borel measure of full support on $[0,1]$. Define
\[\mathcal{E}(f,g)= \int_0^1 f'(x) g'(x) dx,\qquad \forall f,g\in\mathcal{F},\]
where $\mathcal{F}=\{f\in C([0,1]):\: f\mbox{ is absolutely continuous and }f'\in L^2(dx)\}$. Then $(\mathcal{E},\mathcal{F})$ is the resistance form associated with $([0,1],R)$. Moreover, $(\mathcal{E},\mathcal{F})$ is a regular Dirichlet form on $L^2(\mu)$. Integrating by parts yields
\[\mathcal{E}(f,g)= -\int_0^1 (\Delta f)(x) g(x) \mu(dx),\qquad \forall f\in \mathcal{D}(\Delta),\:g\in\mathcal{F},\]
where $\Delta f=\frac{d}{d\mu}\frac{df}{dx}$, and $\mathcal{D}(\Delta)$ contains those $f$ such that: $f'$ exists and $df'$ is absolutely continuous with respect to $\mu$, $\Delta f\in L^2(\mu)$, and $f'(0)=f'(1)=0$. From this, we see that if $\mu(dx)=dx$, then the Markov process naturally associated with $\Delta$ is reflected Brownian motion on $[0,1]$. (For more general $\mu$, the relevant process is simply a time-change of Brownian motion according to Revuz measure $\mu$.) Taking $R(x,y)=|x-y|^{\alpha-1}$ for $\alpha\in(1,2]$, we can also obtain $\alpha$-stable processes in this way (see \cite[Chapter 16]{Kig}).

\section{Scaling limit result}

In this section, we will present a simplified version of the result established in \cite{Croyres}, the aim of which was to establish scaling limits of stochastic processes associated with resistance forms. In the full result, a non-explosion condition was provided to extend from the case of compact spaces that we consider here. Moreover, the result was also adapted to random spaces, and incorporated spatial embeddings. In \cite{CHKtime} a similar result was proved under more restrictive volume growth conditions, which were applied to further deduce a convergence statement regarding the local times of the processes in question.

To introduce the result precisely, let us fix the framework. In particular, we write $\mathbb{F}_c$ for the collection of quadruples of the form $(F,R,\mu,\rho)$, where: $F$ is a non-empty set; $R$ is a resistance metric on $F$ such that $(F,R)$ is compact; $\mu$ is a locally finite Borel regular measure of full support on $(F,R)$; and $\rho$ is a marked point in $F$. We recall that saying a sequence of such spaces converges in the (marked) Gromov-Hausdorff-Prohorov topology to some element of $\mathbb{F}_c$ if all the spaces can be isometrically embedded into a common metric space $(M,d_M)$ in such a way that: the embedded sets converge with respect to the Hausdorff distance, the embedded measures converge weakly, and the embedded marked points converge. The following result establishes that, if such convergence occurs, then we also obtain convergence of stochastic processes.

\begin{theorem}[{cf.\ \cite[Theorem 1.2]{Croyres}}]\label{mainres} Suppose that the sequence $(F_n,R_n,\mu_n,\rho_n)_{n\geq 1}$ in $\mathbb{F}_c$ satisfies
\begin{equation}\label{ghp}
\left(F_n,R_n,\mu_n,\rho_n\right)\rightarrow \left(F,R,\mu,\rho\right)
\end{equation}
in the (marked) Gromov-Hausdorff-Prohorov topology for some $(F,R,\mu,\rho)\in \mathbb{F}_c$. It is then possible to isometrically embed $(F_n,R_n)_{n\geq 1}$ and $(F,R)$ into a common metric space
$(M,d_M)$ in such a way that
\[P^n_{\rho_n}\left(\left(X^n_t\right)_{t\geq 0}\in\cdot\right)\rightarrow P_{\rho}\left(\left(X_t\right)_{t\geq 0}\in\cdot\right)\]
weakly as probability measures on $D(\mathbb{R}_+,M)$ (that is, the space of cadlag processes on $M$, equipped with the usual Skorohod $J_1$-topology), where  $((X^n_t)_{t\geq 0},(P^n_x)_{x\in F_n})$ is the Markov process corresponding to $(F_n,R_n,\mu_n,\rho_n)$, and $((X_t)_{t\geq 0},(P_x)_{x\in F})$ is the Markov process corresponding to $(F,R,\mu,\rho)$.
\end{theorem}

Of course, given the correspondence between measured resistance metric spaces and stochastic processes, as described in Section \ref{ressec}, one might intuitively expect that Gromov-Hausdorff-Prohorov convergence will give us all the information we need to obtain process convergence. To turn this expectation into a proof we use the fact that, for a process associated with a resistance metric, we have an explicit formula for its resolvent kernel. (This starting point was influenced by the one used as the basis of the corresponding argument for trees in \cite{ALWtree}.) In particular, for $(F,R,\mu,\rho)\in \mathbb{F}_c$, let
\[G_xf(y)=E_y\int_0^{\sigma_x}f(X_s)ds\]
be the resolvent of $X$ killed on hitting $x$, where we have written $\sigma_x$ for the hitting time of $x$. (NB.\ Processes associated with resistance forms hit points; the above expression is well-defined and finite.) We then have that
\[G_xf(y)=\int_Fg_x(y,z)f(z)\mu(dz),\]
where the resolvent kernel is given by
\[g_x(y,z)=\frac{R(x,y)+R(x,z)-R(y,z)}{2}.\]
(See \cite[Theorem 4.3]{Kig}.) In view of this expression, the metric measure convergence at (\ref{ghp}) readily gives convergence of resolvents. Relatively standard arguments subsequently yield semigroup convergence, which in turn gives convergence of finite dimensional distributions.

To get from convergence of finite dimensional distributions to convergence in $D(\mathbb{R}_+,M)$, it remains to check tightness of the processes. For this, we again appeal to an explicit expression for a resolvent in terms of resistance. In particular, we have for any closed set $A$ that
\[g_A(y,z)=\frac{R(y,A)+R(z,A)-R_A(y,z)}{2},\]
where $g_A$ is the resolvent kernel for the process $X$ killed on hitting the set $A$, and $R_A(y,z)$ is the resistance from $y$ to $z$ when the set $A$ is `shorted', i.e.\ in defining the resistance between $y$ and $z$ similarly to (\ref{reschar}), we consider only functions that are constant on $A$. (Again, see \cite[Theorem 4.3]{Kig}.) From this, and using that $X$ admits local times $(L_t(x))_{x\in F,t\geq 0}$ (see \cite[Lemma 2.4]{CHKtime}) that satisfy $E_yL_{\sigma_A}(z)=g_A(y,z)$, where $\sigma_A$ is the hitting time of $A$, one can establish via Markov's inequality a general exit time estimate of the form:
\[\sup_{x\in F}P_x\left(\sup_{s\leq t}R(x,X_s)\geq \varepsilon\right)\leq \frac{32 N(F,\varepsilon/4)}{\varepsilon}\left(\delta+\frac{t}{\inf_{x\in F}\mu(B_R(x,\delta))}\right),\]
where $B_R(x,\delta):=\{y\in F:\:R(x,y)<\delta\}$, and $N(F,\varepsilon)$ is the minimal size of an $\varepsilon$ cover of $F$ (see \cite[Lemma 4.3]{Croyres}). The exact form of this expression is not especially important. Rather, the crucial factor is the straightforward dependence on simply metric-measure quantities. As a consequence, we find that if (\ref{ghp}) holds, then
\[\lim_{t\rightarrow0}\limsup_{n\rightarrow\infty}\sup_{x\in F_n}P^n_x\left(\sup_{s\leq t}R_n(x,X^n_s)\geq \varepsilon\right)=0.\]
Tightness of the sequence $X^n$ is then a standard application of Aldous' tightness criterion  \cite[Theorem 16.10 and 16.11]{Kall}.

\section{Applications}\label{applications}

To complete this overview, we present several examples to which the resistance metric framework is particularly useful. For further details, examples and references, see \cite{Croyres, CHKtime}.

\subsection{Trees}

\begin{figure}[ht]
  \begin{center}
\scalebox{.5}{\includegraphics{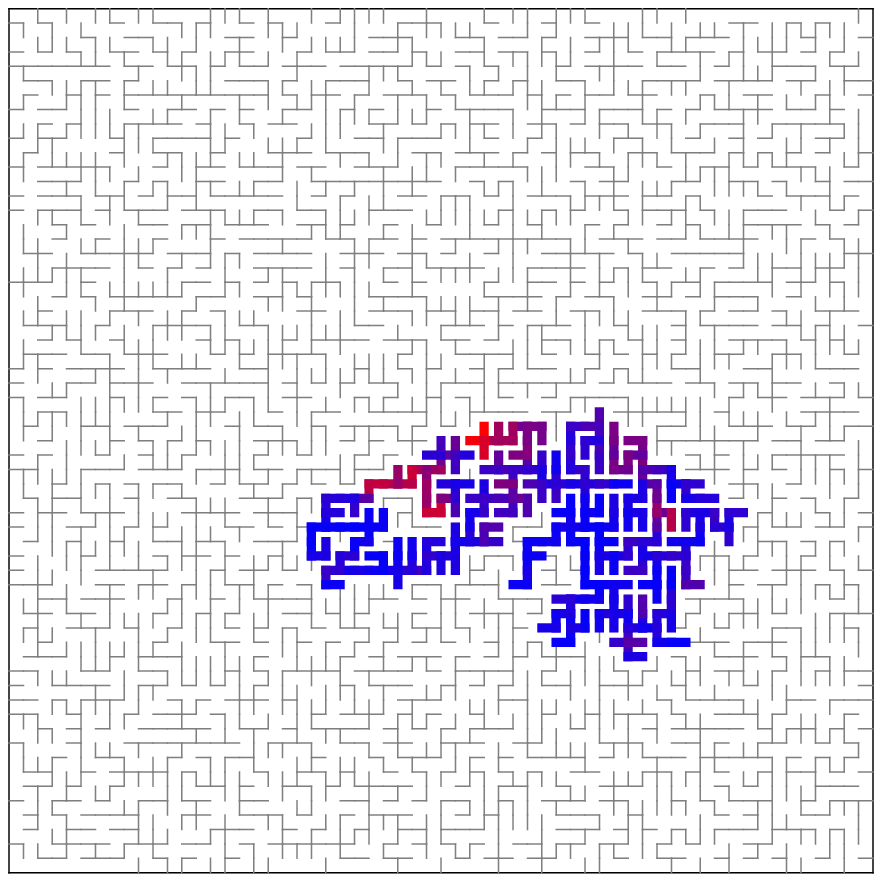} \includegraphics{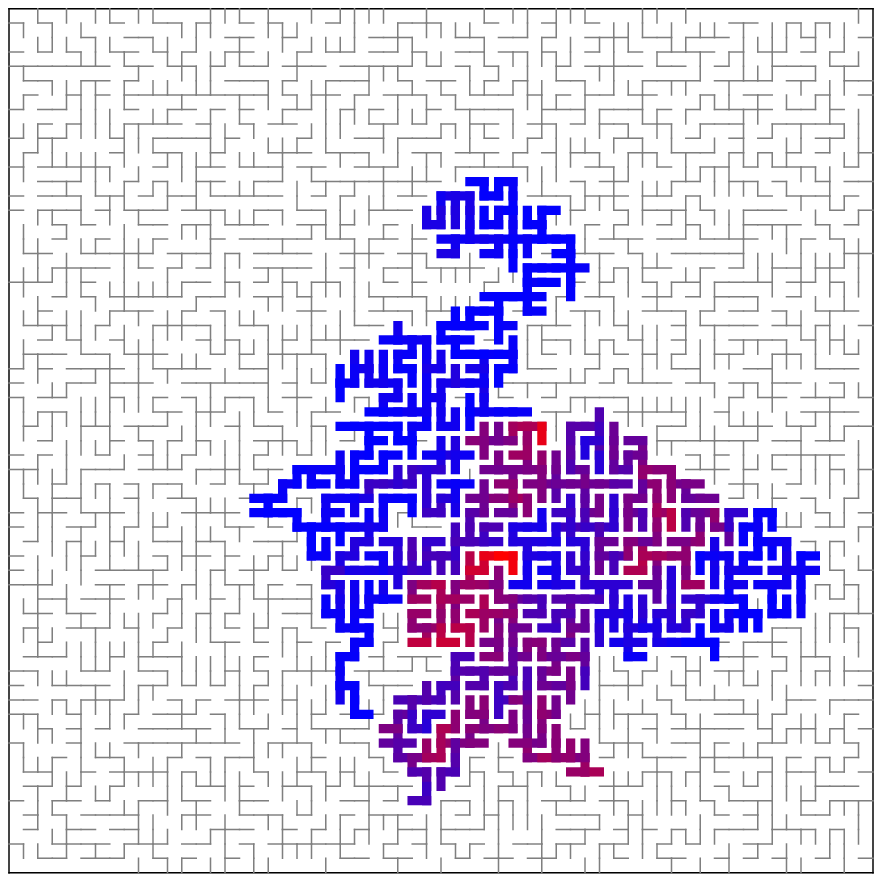}}
   \end{center}
   \caption{The range of a realisation of the simple random walk
on uniform spanning tree on a $60\times60$ box (with wired boundary conditions), shown after 5,000 and 50,000 steps. From most to least crossed edges, colours blend from red to blue. Picture: Sunil Chhita.}\label{ustfig}
\end{figure}

Consider a sequence of graph trees $(T_n)_{n\geq 1}$, where $T_n$ has vertex set $V(T_n)$, shortest path graph distance $R_n$ (noting that this is a resistance metric), counting measure on the vertices $\mu_n$ (placing mass one on each vertex), and root $\rho_n$. Suppose that for some null sequences $(a_n)_{n\geq 1},(b_n)_{n\geq 1}$,
\[\left(V(T_n),a_nR_n,b_n\mu_n,\rho_n\right)\rightarrow\left(\mathcal{T},R,\mu,\rho\right),\]
for some limit in $\mathbb{F}_c$. (NB.\ Under the assumptions stated, $(\mathcal{T},R)$ is a so-called `real tree', which is a natural metric space analogue of a graph tree.) It then holds that
\[\left(X^{T_n}_{ta_nb_n}\right)_{t\geq 0}\rightarrow \left(X^\mathcal{T}_t\right)_{t\geq 0},\]
where $X^{T_n}$ is the random walk associated with $(V(T_n),R_n,\mu_n,\rho_n)$. (Here and in the examples below, for the statement we suppose the state spaces are suitably isometrically embedded into a common metric space.) In particular, distributional versions of the result hold for:
\begin{itemize}
  \item Critical Galton-Watson trees with finite variance conditioned on size, $a_n=n^{1/2}$, $b_n=n$. Versions of the result for infinite variance Galton-Watson trees also hold, see \cite{Croyrogt}.
   \item The (non-lattice) branching random walk, where the underlying tree is a critical Galton-Watson tree with exponential tails for the offspring distribution, and the steps have a centred, continuous distribution with fourth order polynomial tail decay \cite{CroyHd}.
  \item $\Lambda$-coalescent measure trees \cite[Section 7.5]{ALWtree}.
  \item The uniform spanning tree in two dimensions, $a_n=n^{5/4}$, $b_n=n^2$ (this was proved subsequentially in \cite{BCK}, and the full convergence follows from \cite{HoldenSun}). See Figure \ref{ustfig}.
\end{itemize}

\subsection{Conjecture for critical percolation}\label{IICsec}

One model for which an appealing conjecture can be made is the incipient infinite cluster of bond percolation on $\mathbb{Z}^d$ in high dimensions, that is, when $d>6$. In particular, it is expected that this model satisfies the same scaling properties as branching random walk, and thus one might anticipate that if $\mathrm{IIC}$ is the incipient infinite cluster (see \cite{HJ} for a construction), $R_{\mathrm{IIC}}$ is the resistance metric on this (when individual edges have unit resistance) and $\mu_{\mathrm{IIC}}$ is the counting measure on $\mathrm{IIC}$, then the rescaled sequence
\[\left(\mathrm{IIC},n^{-2}R_{\mathrm{IIC}},n^{-4}\mu_{\mathrm{IIC}},0\right),\]
satisfies a locally compact, distributional version of (\ref{ghp}), with limit being (an unbounded version of) the continuum random tree, and so the associated random walks converge to Brownian motion on the latter space, cf.\ the conjecture of \cite{CroyHd}. See also the recent work of \cite{BCF} regarding the lattice branching random walk.

\subsection{Critical random graph}

One critical percolation model that can already be tackled with Theorem \ref{mainres} is that on the complete graph. In particular, let $G(n,1/n)$ be the Erd\H{o}s-R\'{e}nyi random graph at criticality, which is obtained by running bond percolation with edge retention probability $1/n$ on the complete graph with $n$ vertices. For the largest connected component $\mathcal{C}_1^n$ of $G(n,1/n)$, it can be checked that
\[\left(\mathcal{C}_1^n,n^{-1/3}R_n,n^{-2/3}\mu_n,\rho_n\right)\rightarrow\left(F,R,\mu,\rho\right),\]
where the limiting space can be described explicitly, cf.\ \cite{ABG}. Hence, as originally proved in \cite{Croycrg},
\[\left(X^n_{tn}\right)_{t\geq 0}\rightarrow \left(X_t\right)_{t\geq 0}.\]

\subsection{Heavy-tailed random conductance model on fractals}\label{trapsec}

Finally, consider the Sierpinski gasket graphs shown in Figure \ref{sgpic}. Suppose that we equip the edges of these graphs with random, i.i.d.\ edge conductances that satisfy
\[P(c(x,y)\geq u)=u^{-\alpha}\]
for some $\alpha\in (0,1)$. One can then check that resistance homogenises, in the sense that, almost-surely,
\[\left(V_n,(3/5)^nR_n\right)\rightarrow \left(F,R\right),\]
in the Gromov-Hausdorff topology, where: $V_n$ is the vertex set of the $n$th level graph, $R_n$ is the effective resistance associated with the random conductances, $F$ is the Sierpinski gasket, and (up to a deterministic constant) $R$ is the effective resistance on the Sierpinski gasket introduced above, see \cite{CHKtime}.

Recall from Section \ref{rwsec} that the VSRW associated with $(V_n,R_n)$, which has transition rates $\lambda_{xy}= c(x,y)$, is the process corresponding to $(V_n,R_n,\mu_n)$, where $\mu_n(\{x\})=1$. Since $3^{-n}\mu_n\rightarrow \mu$, where $\mu$ is $(\ln3/\ln 2)$-dimensional Hausdorff measure on Sierpinski gasket, it follows that the VSRW $X^n$ converges to the standard Brownian motion on the gasket, $X$ say:
\[\left(X^n_{t5^n}\right)_{t\geq 0}\rightarrow \left(X_t\right)_{t\geq 0}.\]

On the other hand, the associated CSRW has transition rates $\lambda_{xy}= c(x,y)/c(x)$, where $c(x):=\sum_{y:y\sim x}c(x,y)$. This corresponds to the space $(V_n,R_n,\nu_n)$, where $\nu_n(\{x\})=c(x)$. Similarly to the convergence of i.i.d.\ sums to $\alpha$-stable subordinators, it further holds that
\[\nu_n:={3}^{-n/\alpha}\sum_{x\in V_n}c(x)\delta_{x}\rightarrow \nu=\sum_i v_i\delta_{x_i},\]
in distribution, where $\{(v_i,x_i)\}$ is a Poisson point process on $(0,\infty)\times F$ with intensity $c v^{-1-\alpha}dv\mu(dx)$ (for some deterministic constant $c$). Hence the CSRW, $Y^n$ say, (and indeed its jump chain, which is simply the discrete time simple random walk amongst the same conductances) converges:
\[\left(Y^{n}_{t (5/3)^n3^{n/\alpha}}\right)_{t\geq 0}\rightarrow \left(Y_t\right)_{t\geq 0},\]
where the limiting process $Y$ is the so-called \emph{Fontes-Isopi-Newman (FIN)} diffusion on the limiting fractal, which is the time-change of the Brownian motion $X$ according to Revuz measure $\nu$. The process $Y$ spends positive time at atoms of $\nu$, which demonstrates the persistence of the trapping on edges of high conductance of the CSRW in the limit (a phenomenon which the result of the previous paragraph shows does not occur for the related VSRW). For further details of this example, see \cite{CHKtime}.

\providecommand{\bysame}{\leavevmode\hbox to3em{\hrulefill}\thinspace}
\providecommand{\MR}{\relax\ifhmode\unskip\space\fi MR }
\providecommand{\MRhref}[2]{%
  \href{http://www.ams.org/mathscinet-getitem?mr=#1}{#2}
}
\providecommand{\href}[2]{#2}

\end{document}